\documentclass[12pt,a4paper,reqno]{amsart} %Fr o m detta
\usepackage{amssymb}
\usepackage{latexsym}
\usepackage{amsmath}
\usepackage{mathrsfs}
\usepackage{euscript}
\usepackage[X2,T1]{fontenc}
\usepackage{cite}
\usepackage{calc}                   %Detta beh?vs f?r att
\newcommand{\scal}[2]{\langle #1,#2\rangle}
\newcommand{\rr}[1]{\mathbf R^{#1}}

\newcommand{\zz}[1]{\mathbf Z^{#1}}

\newcommand{\nn}[1]{\mathbf N^{#1}}
\newcommand{\nm}[2]{\Vert #1\Vert _{#2}}
\newcommand{\Nm}[2]{\left \Vert #1\right \Vert _{#2}}

\newcommand{\sets}[2]{\{ \, #1\, ;\, #2\, \} }

\newcommand{\fy}{\varphi}
\newcommand{\cdo}{\, \cdot \, }

\newcommand{\vrum}{\vspace{0.1cm}}

\newcommand{\mascB}{\mathscr B}

\newcommand{\mascP}{\mathscr P}

\newcommand{\mabfm}{\boldsymbol m}
\newcommand{\mabfp}{{\boldsymbol p}}

\newcommand{\mabfr}{\boldsymbol r}

\newcommand{\mabfz}{\boldsymbol z}
\newcommand{\mabffy}{\boldsymbol \varphi}

\numberwithin{equation}{section}          %Detta g?r att man f?r
\newtheorem{thm}{Theorem}
\numberwithin{thm}{section}

\newcommand{\rubrik}{}

\theoremstyle{definition}

\newtheorem{defn}[thm]{Definition}

\theoremstyle{remark}
\newtheorem{rem}[thm]{Remark}              %T o m hit  bara allm
\title{Semi-continuous convolutions on weakly periodic
Lebesgue spaces}

\author{Joachim Toft}

\address{Department of Mathematics,
Linn{\ae}us University, V{\"a}xj{\"o}, Sweden}

\email{joachim.toft@lnu.se}

\begin{document}

\par

\begin{abstract}
We deduce mixed quasi-norm estimates of
Lebesgue types on semi-continuous convolutions
between sequences and functions which may be
periodic or possess a weaker form of periodicity
in certain directions. In these directions, the
Lebesgue quasi-norms are applied on the period
instead of the whole axes.
\end{abstract}

\maketitle

%%%%%%%%%%%%%%%%%%%%%%
\section{Introduction}\label{sec0}
%%%%%%%%%%%%%%%%%%%%%%

\par

Continuous, discrete and semi-continuous convolutions appear
naturally when searching for estimates between short-time Fourier
transforms with
different window functions. By straight-forward application
of Fourier's inversion formula, the short-time Fourier
transform $V_\phi f$ of the function or (ultra-)distribution
$f$ with window function $\phi$ is linked to $V_{\phi _0}f$
by
\begin{equation}\label{Eq:STFTConvRel}
|V_\phi f| \lesssim |V_\phi \phi _0|*|V_{\phi _0}f|
\end{equation}
(cf. e.{\,}g. \cite[Chapter 11]{Gc2}). Here $*$ denotes the
usual (continuous) convolution and it is assumed that the
window functions $\phi$ and $\phi _0$ are fixed and belongs
to suitable classes (see \cite{Ho1,Gc2} and
Section \ref{sec1} for notations).

\par

Modulation spaces appear by imposing norm or quasi-norm
estimates on
the short-time Fourier transforms of (ultra-)distributions
in Fourier-invariant spaces. In most situations these
(quasi-)norms are mixed norms of (weighted) Lebesgue types.
More precisely, let $\mascB$ be a mixed quasi-Banach space
of Lebesgue type with functions defined on the phase space,
and let $\omega$ be a moderate weight. Then the modulation
space $M(\omega ,\mascB)$ consists of all ultra-distributions
$f$ such that
\begin{equation}\label{Eq:ModNorm}
\nm f{M(\omega ,\mascB)}\equiv \nm {V_\phi f \cdot \omega}{\mascB}
\end{equation}
is finite.
% The modulation spaces were introduced by Feichtinger
% in \cite{Fe2}, and extended in several ways in e.{\,}g. 
% \cite{Fe3,GH1,GH2}.

\par

If $\mascB$ is a Banach space of mixed Lebesgue type, then the
inequality \eqref{Eq:STFTConvRel} can be used to deduce:
\begin{enumerate}
\item that $M(\omega ,\mascB )$ is invariant of the choice of
window function $\phi$ in \eqref{Eq:ModNorm}, and that different
$\phi$ give rise to equivalent norms.

\vrum

\item that $M(\omega ,\mascB)$ increases with the Lebesgue
exponents.

\vrum

\item that $M(\omega ,\mascB)$ is complete.
\end{enumerate}

\par

Essential parts of these basic properties for modulation spaces
were established in the pioneering paper \cite{Fe2}, but some tracks
goes back to \cite{Fe1}. The theory has thereafter been developed
in different ways, see e.{\,}g. \cite{Fe3,FG1,FG2,Gc2}.

\par

A more complicated situation appear when some of the Lebesgue parameters
for $\mascB$ above are strictly smaller than one, since $\mascB$
is then merely a quasi-Banach space, but not a Banach space, since
only a weaker form of the triangle inequality hold true. In such
situation, $\mascB$ even fails to be a local convex topological vector
space, and the analysis based on \eqref{Eq:STFTConvRel} to reach (1)--(3)
in their full strength above seems not work. (Some partial properties
can be achieved if for example it is required that the Fourier
transform of $\phi$ and $\phi _0$ should be compactly supported, see
e.{\,g} \cite{WaHu}.)

\par

In \cite{GaSa}, the more discrete approach is used to handle this
situation, where a Gabor expansion of $\phi$ with $\phi _0$ as
Gabor window leads to that $|V_\phi f|$ can be estimated by
\begin{equation}\label{Eq:STFTConvRe2}
|V_\phi f| \lesssim a*_{[E]}|V_{\phi _0}f|,
\end{equation}
for some non-negative sequence $a$ with enough rapid decay
towards zero at infinity.
Here $*_{[E]}$ denotes the semi-continuous convolution
$$
a*_{[E]}F \equiv \sum _{j\in \Lambda _E}F(\cdo -j)a(j)
$$
with respect to the basis $E$, between functions
$F$ and sequences $a$, and $\Lambda _E$ is the lattice spanned by
$E$. It follows that $*_{[E]}$ is similar to discrete convolutions.

\par

For the discrete convolution $*$ both the classical
Young's inequality 
\begin{alignat}{3}
\nm {a*b}{\ell ^p_0} &\le \nm a{\ell ^{p_1}}\nm b{\ell ^{p_2}},&
\quad
\frac 1{p_1}+\frac 1{p_2} &= 1+\frac 1{p_0},&\ p_j&\in [1,\infty ],
\intertext{as well as}
\nm {a*b}{\ell ^p} &\le \nm a{\ell ^p}\nm b{\ell ^r},&
\quad r&\le \min (1,p),& \ \ p,r&\in (0,\infty ],
\end{alignat}
hold true, and it is proved in \cite{GaSa} and
extended in \cite{To15} that similar facts hold true
for semi-continuous convolutions. In the end
the following restatement of \cite[Proposition 2.1]{To15} is deduced.
The result also extends \cite[Lemma 2.6]{GaSa}.

\par

\begin{thm}\label{Thm:SemiContConvEst0}
Let $E$ be an ordered basis
of $\rr d$,
$\omega ,v\in \mascP _E(\rr d)$ be
such that $\omega$ is $v$-moderate, and let
$\mabfp ,\mabfr \in (0,\infty ]^{d}$ be such that
$$
r_k\le \min _{m\le k}(1,p_m).
$$
Also let 
$f$ be measurable.
Then the map $(a,f)\mapsto a*_{[E]}f$ from $\ell _0(\Lambda _E)\times
\Sigma _1(\rr d)$ to
$L^{\mabfp}_{E,(\omega )}(\rr d)$ extends uniquely to a
linear and continuous map from $\ell ^{\mabfr}_{E ,(v)}(\Lambda _E)
\times L^{\mabfp}_{E,(\omega )}(\rr d)$ to
$L^{\mabfp}_{E,(\omega )}(\rr d)$, and
\begin{equation}\label{convest0}
\nm {a*_{[E]}f}{L^{\mabfp }_{E,(\omega )}}\lesssim
\nm {a}{\ell ^{\mabfr}_{E ,(v)}}
\nm {f}{L^{\mabfp }_{E,(\omega )}}.
\end{equation}
\end{thm}

\par

In \cite{GaSa}, \eqref{Eq:STFTConvRe2} in combination with
\cite[??]{GaSa} is used to show that (1)--(3) still hold
when $\mascB =L^{p,q}$ and $\omega$ is a moderate weight of
polynomial type. In \cite{To15}, \eqref{Eq:STFTConvRe2} in
combination with Theorem \ref{Thm:SemiContConvEst0} are used
to show (1)--(3) for an even broader class of mixed
Lebesgue spaces $\mascB$ and weight functions $\omega$.

\par

The aim of the paper is to extend Theorem
\ref{Thm:SemiContConvEst0}, so that $f$ in some directions
(variables) is allowed to be periodic, or a weaker form of
periodicity, called \emph{echo-periodic functions}. Such functions
appear for example when applying the short-time Fourier transform on
periodic or quasi-periodic functions. In fact, if
$f$ is $E$-periodic, then $x\mapsto |V_\phi f(x,\xi )|$
is $E$-periodic for every $\xi$. A function or distribution
$F(x,\xi )$ is called quasi-periodic of order $\rho >0$, if
\begin{equation*}
\begin{alignedat}{2}
F(x+\rho k,\xi ) &= e^{2\pi i\rho \scal k\xi} F(x,\xi ),& \quad k &\in \zz d,
\\[1ex]
F(x,\xi +\kappa /\rho ) &= F(x,\xi ),& \quad \kappa &\in \zz d.
\end{alignedat}
\end{equation*}
and by straight-forward computations it follows that
\begin{alignat}{2}
|(V_\Phi F)(x+\rho k,\xi ,\eta ,y)| &=
|(V_\Phi F)(x,\xi ,\eta ,y-2\pi k)|,
& \quad k &\in \zz d,\label{Eq:PerTransfer}
\\[1ex]
|(V_\Phi F)(x,\xi +\kappa /\rho ,\eta ,y)| &=
|(V_\Phi F)(x,\xi ,\eta ,y)|,
&\quad \kappa &\in \zz d,\notag
\end{alignat}
for such $F$.

\par

It is expected that
the achieved extensions will be useful when performing
local investigations of short-time Fourier transforms
of periodic and quasi-periodic functions, e.{\,}g. in
\cite{To18}.

\par

%%%%%%%%%%%%%%%%%%%%%%
\section{Preliminaries}\label{sec1}
%%%%%%%%%%%%%%%%%%%%%%

\par

In this section we recall some basic facts and introduce
some notations. In the first part we recall the notion
of weight functions. Thereafter we discuss mixed quasi-norm
spaces of Lebesgue types. Finally we consider periodic
functions and distributions, and introduce the notion of
echo-periodic functions, which is a weaker form of periodicity
which at the same time also include the notion of
quasi-periodicity.

\par

\subsection{Weight functions}\label{subsec1.1}

\par

A \emph{weight} on $\rr d$ is a positive function $\omega
\in  L^\infty _{loc}(\rr d)$ such that $1/\omega \in  L^\infty _{loc}(\rr d)$.
A usual condition on $\omega$ is that it should be \emph{moderate},
or \emph{$v$-moderate} for some positive function $v \in
 L^\infty _{loc}(\rr d)$. This means that
\begin{equation}\label{moderate}
\omega (x+y) \lesssim \omega (x)v(y),\qquad x,y\in \rr d.
\end{equation}
We note that \eqref{moderate} implies that $\omega$ fulfills
the estimates
\begin{equation}\label{moderateconseq}
v(-x)^{-1}\lesssim \omega (x)\lesssim v(x),\quad x\in \rr d.
\end{equation}
We let $\mascP _E(\rr d)$ be the set of all moderate weights on $\rr d$.

\par

It can be proved that if $\omega \in \mascP _E(\rr d)$, then
$\omega$ is $v$-moderate for some $v(x) = e^{r|x|}$, provided the
positive constant $r$ is large enough (cf. \cite{Gc2.5}). In particular,
\eqref{moderateconseq} shows that for any $\omega \in \mascP
_E(\rr d)$, there is a constant $r>0$ such that
$$
e^{-r|x|}\lesssim \omega (x)\lesssim e^{r|x|},\quad x\in \rr d.
$$

\par

We say that $v$ is
\emph{submultiplicative} if $v$ is even and \eqref{moderate}
holds with $\omega =v$. In the sequel, $v$ and $v_j$ for
$j\ge 0$, always stand for submultiplicative weights if
nothing else is stated.

\par

\subsection{Spaces of mixed quasi-norm spaces of
Lebesgue types}\label{subsec1.2}

Our discussions on periodicity are done in terms of suitable bases.

\par

\begin{defn}\label{Def:OrdBasis}
Let $E$ be an (ordered) basis $e_1,\dots,e_d$ to
$\rr {d}$. Then
\begin{align*}
\Lambda _E
&=
\sets{n_1e_1+\cdots+n_de_d}{(n_1,\dots,n_d)\in \zz d}
\end{align*}
is the corresponding lattices.
\end{defn}

% \begin{defn}\label{Def:OrdBasis}
% Let $E$ be an (ordered) basis $e_1,\dots,e_d$ to
% $\rr {d}$. Then $E'$ denotes the basis of $e_1',\dots,e_d'$ 
% in $\rr {d}$ which satisfies
% $$
% \scal {e_j} {e'_k}=2\pi \delta_{jk}
% \quad \text{for every}\quad
% j,k =1,\dots, d,
% $$
% and
% %%
% \begin{align*}
% \Lambda _E
% &=
% \sets{n_1e_1+\cdots+n_de_d}{(n_1,\dots,n_d)\in \zz d},
% \intertext{and}
% \Lambda'_E
% &=
% \Lambda_{E'}=\sets{\nu _1e'_1+\cdots+\nu _de'_d}{(\nu _1,\dots,\nu _d)
% \in \zz d},
% \end{align*}
% %%
% are the corresponding lattices. The sets $E'$ and $\Lambda '_E$ are
% called the dual basis and dual lattice of $E$ and $\Lambda _E$, 
% respectively.
% \end{defn}

\par

Evidently, if $E$ is the same as in Definition \ref{Def:OrdBasis},
then there is a matrix $T_E$ with $E$ as the image of
the standard basis in $\rr d$. Then $E'$ is the image of the standard basis
under the map $T_{E'}= 2\pi(T^{-1}_E)^t$.

\par

\begin{defn}\label{Def:DiscLebSpaces}
Let $E$ be a basis of $\rr d$, $\kappa (E)$ be the
parallelepiped spanned by $E$, $\omega \in \mascP _E(\rr d)$
$\mabfp =(p_1,\dots ,p_d)\in (0,\infty ]^{d}$ and $r=\min (1,\mabfp )$.
If  $f\in L^r_{loc}(\rr d)$, then
$$
\nm f{L^{\mabfp }_{E,(\omega )}}\equiv
\nm {g_{d-1}}{L^{p_{d}}(\mathbf R)}
$$
where  $g_k(\mabfz _k)$, $z_k\in \rr {d-k}$,
$k=0,\dots ,d-1$, are inductively defined as
\begin{align*}
g_0(x_1,\dots ,x_{d})
&\equiv
|f(x_1e_1+\cdots +x_{d}e_d)\omega (x_1e_1+\cdots +x_{d}e_d)|,
\\[1ex]
\intertext{and}
g_k(\mabfz _k) &\equiv
\nm {g_{k-1}(\cdo ,\mabfz _k)}{L^{p_k}(\mathbf R)},
\quad k=1,\dots ,d-1.
\end{align*}
\begin{enumerate}
\item If $\Omega \subseteq \rr d$ is measurable,
then $L^{\mabfp }_{E,(\omega )}(\Omega )$ consists
of all $f\in L^r_{loc}(\Omega )$ with finite quasi-norm
$$
\nm f{L^{\mabfp}_{E,(\omega )}(\Omega )}
\equiv
\nm {f_\Omega }{L^{\mabfp}_{E,(\omega )}(\rr d)},
\qquad
f_\Omega (x)
\equiv
\begin{cases}
f(x), &\text{when}\ x\in \Omega
\\[1ex]
0, &\text{when}\ x\notin \Omega .
\end{cases}
$$
The space $L^{\mabfp }_{E,(\omega )}(\Omega )$ is called 
\emph{$E$-split Lebesgue space (with respect to $\omega$, $\mabfp$,
$\Omega$ and $E$)};

\vrum

\item If $\Lambda \subseteq \rr d$ is a lattice such that
$\Lambda _E\subseteq \Lambda$, then the quasi-Banach space
$\ell ^{\mabfp } _{E ,(\omega )}(\Lambda )$ consists of all
$a\in \ell _0'(\Lambda )$ such that
$$
\nm a{\ell ^{\mabfp }_{E,(\omega )}(\Lambda )}
\equiv
\Nm {\sum _{j\in \Lambda}a(j)\chi _{j+\kappa (E)}}
{L^{\mabfp }_{E,(\omega )}(\rr d)}
$$
is finite. The space $\ell ^{\mabfp }_{E,(\omega )} \equiv \ell ^{\mabfp }
_{E,(\omega )}(\Lambda _E)$ is called the
\emph{discrete version of $L^{\mabfp }_{E,(\omega )}(\rr d)$}.
\end{enumerate}
\end{defn}

\par

Evidently, $L^{\mabfp}_{E,(\omega )} (\Omega )$ and
$\ell ^{\mabfp}_{E,(\omega )} (\Lambda )$
in Definition \ref{Def:DiscLebSpaces} are quasi-Banach spaces of order
$\min (\mabfp ,1)$. We set
$$
L^{\mabfp}_{E} = L^{\mabfp}_{E,(\omega )}
\quad \text{and}\quad
\ell ^{\mabfp}_{E} = \ell ^{\mabfp}_{E,(\omega )}
$$
when $\omega =1$, and if $\mabfp = (p,\dots ,p)$ for some
$p\in (0,\infty ]$, then
\begin{alignat*}{5}
L^{p}_{E,(\omega )} &=  L^{\mabfp}_{E,(\omega )},
&\quad
L^{p}_{E} &= L^{\mabfp}_{E},
&\quad
\ell ^{p}_{E,(\omega )} &= \ell ^{\mabfp}_{E,(\omega )}
&\quad &\text{and} &\quad
\ell ^{p}_{E} &= \ell ^{\mabfp}_{E}
\intertext{agree with}
&L^p_{(\omega )}, &\qquad
&L^p, &\qquad
&\ell ^{p}_{(\omega )}
&\quad &\text{and}&\quad
&\ell ^{p},
\end{alignat*}
respectively, with equivalent quasi-norms.

\par

\subsection{Periodic and echo-periodic functions}

\par

We recall that if $E =\{ e_1,\dots ,e_d\}$ is an ordered
basis of $\rr d$, then
the function or distribution $f$ on $\rr d$ is called $E$-periodic, if
$f(\cdo +v)=f$ for every $v\in E$. More generally, if $E_0=\subseteq E$,
then $f$ above is called $E_0$-periodic, if
$f(\cdo +v)=f$ for every $v\in E_0$. We shall consider functions that
possess weaker periodic like conditions, which appear when dealing with
e.{\,}g. quasi-periodic functions and their short-time Fourier transforms.

\par

\begin{defn}\label{Def:PerEcho}
Let $E =\{ e_1,\dots ,e_d\}$ be an ordered basis of $\rr d$,
$E_0\subseteq E$
and let $f$ be a (complex-valued) function on $\rr d$. For every
$k\in \{ 1,\dots ,d\}$, let
$M_{k}$ be the set of all $l\in \{ 1,\dots ,k\}$ such that
$e_l\in E\setminus E_0$.
Then $f$ is called an \emph{echo-periodic function with
respect to $E_0$}, if for every
$e_k\in E_0$, there is a
vector
$$
v_k = \sum _{l\in M_k} v_{k,l}e_l
$$
such that
\begin{equation}\label{Eq:PerEchoDef}
|f(\cdo +e_k)|=|f(\cdo +v_k)|.
\end{equation}
\end{defn}

\par

We notice that in \eqref{Eq:PerTransfer}, relations of
the form \eqref{Eq:PerEchoDef} appears.

\par

\begin{rem}\label{Rem:PerEcho}
Let $E$, $E_0$ and $M_k$ be the same as in Definition \ref{Def:PerEcho},
and let  $f$ be a (complex-valued) function on $\rr d$
such that \eqref{Eq:PerEchoDef} holds true. Also let
\begin{alignat*}{1}
J_k
&=
{
\begin{cases}
\mathbf R,& k\in M_d,
\\[1ex]
[0,1],& k\notin M_d,
\end{cases}
}
\\[1ex]
I_k &=  \sets {xe_k}{x\in J_k},\qquad k\in \{1,\dots ,d\}
\intertext{and}
I &= \sets {x_1e_1+\cdots x_de_d}{x_k\in J_k,\ k=1,\dots ,d}
\\
&\simeq I_1\times \cdots \times I_d.
\end{alignat*}
Then evidently,
$|f(\cdo +ne_k)|=|f(\cdo +nv_k)|$ for every integer $n$.
Hence, if $f$
is measurable and echo-periodic with respect to
$E_0$, and $\mabfp \in (0,\infty ]^d$, then
it follows by straight-forward computations that
$$
\nm {f(\cdo +ne_k)}{L^{\mabfp}_E(I)} = \nm f{L^{\mabfp}_E(I)}
$$
for every integer $n$ and $e_k\in E_0$.
\end{rem}

\par

\begin{defn}\label{Def:PerEchoLebSpaces}
Let $E$, $E_0$ and $I\subseteq \rr d$ be the same as in Remark 
\ref{Rem:PerEcho}, $\omega \in \mascP _E(\rr d)$
and let $\mabfp \in (0,\infty ]^d$. Then
$L^{\mabfp ,E_0}_{E,(\omega )}(\rr d)$ denotes the set of
all complex-valued measurable echo-periodic functions
$f$ with respect to $E_0$ such that
$$
\nm f{L^{\mabfp ,E_0}_{E,(\omega )}}
\equiv
\nm f{L^{\mabfp}_{E,(\omega )}(I)}
$$
is finite.
\end{defn}

\par

In the next section we shall deduce weighted $L^{\mabfp}_E(I)$
estimates of the \emph{semi-discrete convolution}
\begin{equation}\label{EqDistSemContConv}
(a*_{[E]}f)(x) \sum _{j\in \Lambda _E}a(j)f(x - j),
\end{equation}
of the measurable function $f$ on $\rr d$ and
$a \in \ell _0 (\Lambda _E)$,
with respect to the ordered basis $E$.

\par

%%%%%%%%%%%%%%%%%%%%%%%%%%%%%%%%%%%
\section{Weighted Lebesgue estimates
on semi-discrete convolutions}\label{sec2}
%%%%%%%%%%%%%%%%%%%%%%%%%%%%%%%%%%%

\par

In this section we extend Theorem \ref{Thm:SemiContConvEst0}
from the introduction such that $L^{\mabfp}_E(I)$-estimates
of echo-periodic functions are included.

\par

Let $E$, $E_0$, $M_k$ and $J_k$, $k=1,\dots ,d$,
be the same as in Remark \ref{Rem:PerEcho}.
In what follows we let $\Sigma _1^{E_0}(\rr d)$ be the set
of all $E_0$-periodic $f\in C^\infty (\rr d)$ such that if
$$
g(x_1,\dots ,x_d)\equiv f(x_1e_1+\cdots +x_de_d),
$$
then
$$
\sup _{\alpha ,\beta \in \nn d}
\frac {\nm {x^\alpha D^\beta g}{L^\infty (I)}}
{h^{|\alpha +\beta|}\alpha !\beta !}
$$
is finite for every $h>0$. By the assumptions and
basic properties due to \cite{ChuChuKim} it follows that
$\Sigma _1^{E_0}(\rr d)\subseteq
L^{\mabfp ,E_0}_{E,(\omega )}(\rr d)$ for every choice
of $\omega \in \mascP _E(\rr d)$ and $\mabfp \in (0,\infty ]^d$
such that
\begin{equation}\label{Eq:E0WeightCond}
\omega (x)=\omega (x_0)
\qquad \text{when}\qquad
x = \sum _{k=1}^d x_ke_k,
\quad
x_0 = \sum _{k\in M_d} x_ke_k.
\end{equation}

\par

Our extension of Theorem \ref{Thm:SemiContConvEst0}
to include echo-periodic functions is the following,
which is also our main result.

\par

\begin{thm}\label{Thm:SemiContConvEst}
Let $E$ be an ordered basis
of $\rr d$, $E_0\subseteq E$,
$\omega ,v\in \mascP _E(\rr d)$ be
such that $\omega$ is $v$-moderate and satisfy
\eqref{Eq:E0WeightCond}, and let
$\mabfp ,\mabfr \in (0,\infty ]^{d}$ be such that
$$
r_k\le \min _{m\le k}(1,p_m).
$$
Also let 
$f$ be measurable echo-periodic function with respect to
$E_0$, and let
$I\subseteq \rr d$ be as in Remark \ref{Rem:PerEcho}.
Then the map $(a,f)\mapsto a*_{[E]}f$ from $\ell _0(\Lambda _E)\times
\Sigma _1^{E_0}(\rr d)$ to
$L^{\mabfp}_{E,(\omega )}(I)$ extends uniquely to a
linear and continuous map from $\ell ^{\mabfr}_{E ,(v)}(\Lambda _E)
\times L^{\mabfp ,E_0}_{E,(\omega )}(\rr d)$ to
$L^{\mabfp}_{E,(\omega )}(I)$, and
\begin{equation}\label{convest1}
\nm {a*_{[E]}f}{L^{\mabfp }_{E,(\omega )}(I)}\lesssim
\nm {a}{\ell ^{\mabfr}_{E ,(v)}(\Lambda _E)}
\nm {f}{L^{\mabfp }_{E,(\omega )}(I)}.
\end{equation}
\end{thm}

\par

For the proof we recall that
\begin{equation}\label{Eq:OtherMinkowski}
\left ( 
\sum _{j\in I} |b(j)|
\right )^{r}
\le
\sum _{j\in I} |b(j)|^r
\qquad
0<r\le 1,
\end{equation}
for any sequence $b$ and countable set $I$.

\par

\begin{proof}
By letting
\begin{align*}
f_0(x_1,\dots ,x_d)
&=
|f(x_1e_1+\cdots +x_de_d)\omega (x_1e_1+\cdots +x_de_d)|,
\\[1ex]
a_0(l_1,\dots ,l_d)
&=
|a(l_1e_1+\cdots +l_de_d)v(l_1e_1+\cdots +l_de_d)|
\end{align*}
and using the inequality
$$
|a*_{[E]}f \cdot \omega |\lesssim a_v *_{[E]}f_\omega ,
$$
we reduce ourselves to the case when $E$ is the standard
basis, $\omega =v=1$ and $f,a\ge 0$. This implies that we may identify
$I_k$ in Remark \ref{Rem:PerEcho} with $J_k$ for every $k$.

\par

Let
\begin{gather*}
\mabfz _k = (x_{k+1},\dots ,x_d)\in \rr {d-k},\qquad
\mabfm _k = (l_{k+1},\dots ,l_d)\in \zz {d-k}
\intertext{for $k=0,\dots ,d-1$, and let}
f_0 = f,\qquad a_0=a,\qquad g_0= a*_{[E]}f.
\end{gather*}
Then $\mabfz _{k-1}=(x_k,\mabfz _k)$ and
$\mabfm _{k-1}=(l_k,\mabfm _k)$. It follows that
$x_k\in I_k$ when applying the mixed quasi-norms of
Lebesgue types, and that 
\begin{multline}\label{Eq:ConvRef1}
0\le (a*_{[E]}f)(x_1,\dots ,x_d)
\\[1ex]
\le
\sum _{\mabfm _0\in \zz d}
f(x_1-\fy _1(\mabfm _0),\dots ,z_d-\fy _d(\mabfm _{d-1}))a(\mabfm _0),
\end{multline}
for some linear functions $\fy _k$ from
$\rr {d+1-k}$ to $\mathbf R$, which satisfy
\begin{equation}\label{Eq:PhiSeqDef}
\fy _k(\mabfz _{k-1})
=
\begin{cases}
x_k +\psi _k(\mabfz _k),& J_k=\mathbf R,
\\[1ex]
0, & J_k=[0,1],
\end{cases}
\end{equation}
for some linear forms $\psi _k$ on $\rr {d-k}$, $k=1,\dots ,d$.

\par

Define inductively
\begin{align*}
f_k(\mabfz _k) &= \nm {f_{k-1}(\cdo ,\mabfz _k)}{L^{p_k}(J_k)},
\quad
a_k(\mabfm _k) = \nm {a_{k-1}(\cdo ,\mabfm _k)}{\ell ^{r_k}(\mathbf Z)},
\intertext{and}
g_k(\mabfz _k) &= \nm {g_{k-1}(\cdo ,\mabfz _k)}{L^{p_k}(J_k)},\qquad k=1,\dots d.
\end{align*}
Also let 
$$
\mabffy _k(\mabfz _k) = (\fy _{k+1}(\mabfz _{k}),\dots ,\fy _d(\mabfz _{d-1})),
\quad k=0,\dots ,d-1.
$$
Then \eqref{Eq:ConvRef1} is the same as
\begin{equation}\label{Eq:ConvRef2}
0\le (a*_{[E]}f)(x_1,\dots ,x_d)
\le
\sum _{\mabfm _0\in \zz d}
f(\mabfz _0 -\mabffy _0(\mabfm _0))a(\mabfm _0),
\end{equation}

\par

We claim
$$
g_k(\mabfz _k)
\lesssim
\left (
\sum _{\mabfm _k} f_k(x_{k+1}-\fy _{k+1}(\mabfm _k),
\dots ,x_{d}-\fy _{d}(\mabfm _{d-1}))^{p_{0,k}}a_k(\mabfm _k)^{p_{0,k}}
\right )^{\frac 1{p_{0,k}}},
$$
which in view of the links between \eqref{Eq:ConvRef1} and
\eqref{Eq:ConvRef2} is the same as
\begin{equation}\label{Eq:ConvIndEst1}
g_k(\mabfz _k)
\lesssim
\left (
\sum _{\mabfm _k} f_k(\mabfz _{k}-\mabffy _{k}(\mabfm _k))^{p_{0,k}}
a_k(\mabfm _k)^{p_{0,k}}
\right )^{\frac 1{p_{0,k}}}
\end{equation}
when $k=0, \dots ,d$. Here we set $p_{0,0}=1$, and
interprete $f_d$, $a_d$, $g_d$ and the right-hand side of
\eqref{Eq:ConvIndEst1} as $\nm f{L^{\mabfp}_E(I)}$, $\nm a{\ell
^{\mabfp _0}_E(\zz d)}$, $\nm {g_0}{L^{\mabfp}_E(I)}$ and
$\nm f{L^{\mabfp}_E(I)}\nm a{\ell ^{\mabfp _0}_E(\zz d)}$,
respectively. The result then follows by letting $k=d$ in
\eqref{Eq:ConvIndEst1}.

\par

We shall prove \eqref{Eq:ConvIndEst1} by induction. The result
is evidently true when $k=0$. Suppose it is true for $k-1$,
$k\in \{1,\dots ,d-1\}$. We shall consider the cases when
$p_k\ge p_{0,k-1}$ or $p_k\le p_{0,k-1}$, and $J_k=\mathbf R$ or $J_k=[0,1]$
separately, and for conveniency we set $p_{0,k-1}=p$ and $f_{k-1}=h$.

\par

First assume that $p_k\ge p_{0,k-1}$. Then $p_{0,k}=p_{0,k-1}$.
Also suppose $J_k=\mathbf R$. Then
it follows from the induction
hypothesis that
\begin{multline*}
g_k(\mabfz _k)
\\[1ex]
\lesssim
\left (
\int _{-\infty}^\infty
\left (
\sum
h(x_k-\fy _k(\mabfm _{k-1}), \mabfz _k-\mabffy _k(\mabfm _k))
^{p}a_{k-1}(l_k,\mabfm _k)^{p}
\right )^{\frac {p_k}{p}}\, dx_k
\right )^{\frac 1{p_k}},
\end{multline*}
where the sum is taken over all $(l_k,\mabfm _k)\in \mathbf Z\times \zz {d-k}$.
By Minkowski's inequality, the right-hand side can be estimated by
\begin{multline*}
\left (
\sum 
\left (
\int _{-\infty}^\infty
h(x_k-\fy _k(\mabfm _{k-1}), \mabfz _k-\mabffy _k(\mabfm _k))
^{p_{k}}\, dx_k
\right )^{\frac {p}{p_k}}
a_{k-1}(l_k,\mabfm _k)^{p}
\right )^{\frac 1{p}}
\\[1ex]
=
\left (
\sum 
\left (
\int _{-\infty}^\infty
h(x_k, \mabfz _k-\mabffy _k(\mabfm _k))
^{p_{k}}\, dx_k
\right )^{\frac {p}{p_k}}
a_{k-1}(l_k,\mabfm _k)^{p}
\right )^{\frac 1{p}}
\\[1ex]
=
\left (
\sum 
f_{k}(\mabfz _k-\mabffy _k(\mabfm _k)) ^{p}
a_{k-1}(l_k,\mabfm _k)^{p}
\right )^{\frac 1{p}}
\\[1ex]
=
\left (
\sum _{\mabfm _k\in \zz {d-k}}
f_{k}(\mabfz _k-\mabffy _k(\mabfm _k)) ^{p}
\left (
\sum _{l_k\in \mathbf Z}a_{k-1}(l_k,\mabfm _k)^{p}
\right )
\right )^{\frac 1{p}}
\\[1ex]
=
\left (
\sum _{\mabfm _k\in \zz {d-k}}
f_{k}(\mabfz _k-\mabffy _k(\mabfm _k)) ^{p_{0,k}}
a_{k}(\mabfm _k)^{p_{0,k}}
\right )^{\frac 1{p_{0,k}}},
\end{multline*}
and \eqref{Eq:ConvIndEst1} follows in the case $p_k\ge p_{0,k-1}$
and $J_k=\mathbf R$ by combining these estimates.

\par

Next we consider the case when $p_k\ge p_{0,k-1}$
and $J_k=[0,1]$. Then $\fy _k(\mabfm _{k-1})=0$, and
by the induction hypothesis and Minkowski's
inequality we get
\begin{multline*}
g_k(\mabfz _k)
\\[1ex]
\lesssim
\left (
\int _{0}^1
\left (
\sum _{\mabfm _{k-1}}
h(x_k, \mabfz _k-\mabffy _k(\mabfm _k))
^{p}a_{k-1}(l_k,\mabfm _k)^{p}
\right )^{\frac {p_k}{p}}\, dx_k
\right )^{\frac 1{p_k}}
\\[1ex]
\le
\left (
\sum _{\mabfm _{k-1}}
\left (
\int _{0}^1
h(x_k, \mabfz _k-\mabffy _k(\mabfm _k))
^{p_{k}}\, dx_k
\right )^{\frac {p}{p_k}}
a_{k-1}(l_k,\mabfm _k)^{p}
\right )^{\frac 1{p}}
\end{multline*}
\begin{multline*}
=
\left (
\sum _{\mabfm _{k-1}}
f_{k}(\mabfz _k-\mabffy _k(\mabfm _k)) ^{p}
a_{k-1}(l_k,\mabfm _k)^{p}
\right )^{\frac 1{p}}
\\[1ex]
=
\left (
\sum _{\mabfm _k\in \zz {d-k}}
f_{k}(\mabfz _k-\mabffy _k(\mabfm _k)) ^{p_{0,k}}
a_{k}(\mabfm _k)^{p_{0,k}}
\right )^{\frac 1{p_{0,k}}},
\end{multline*}
and \eqref{Eq:ConvIndEst1} follows in the case $p_k\ge p_{0,k-1}$
and $J_k=[0,1]$ as well.

\par

Next assume that $p_k\le p_{0,k-1}$
and $J_k=\mathbf R$. Then
$$
p_k/p_{0,k-1}=p_k/p\le 1
\quad \text{and}\quad
p_{0,k}=p_k,
$$
and \eqref{Eq:OtherMinkowski} gives
\begin{multline*}
g_k(\mabfz _k)
\\[1ex]
\lesssim
\left (
\int _{-\infty}^\infty
\left (
\sum _{\mabfm _{k-1}}
h(x_k-\fy _k(\mabfm _{k-1}), \mabfz _k-\mabffy _k(\mabfm _k))
^{p}a_{k-1}(l_k,\mabfm _k)^{p}
\right )^{\frac {p_k}{p}}\, dx_k
\right )^{\frac 1{p_k}}
\\[1ex]
\lesssim
\left (
\int _{-\infty}^\infty
\sum _{\mabfm _{k-1}}
\left (
h(x_k-\fy _k(\mabfm _{k-1}), \mabfz _k-\mabffy _k(\mabfm _k))
^{p}a_{k-1}(l_k,\mabfm _k)^{p}
\right )^{\frac {p_k}{p}}\, dx_k
\right )^{\frac 1{p_k}}
% \end{multline*}
% %%
% %%
% \begin{multline*}
\\[1ex]
=
\left (
\sum _{\mabfm _{k-1}}
\left (
\int _{-\infty}^\infty
h(x_k-\fy _k(\mabfm _{k-1}), \mabfz _k-\mabffy _k(\mabfm _k))
^{p_k}\, dx_k
\right )
a_{k-1}(l_k,\mabfm _k)^{p_k}
\right )^{\frac 1{p_k}}
\\[1ex]
=
\left (
\sum _{\mabfm _{k-1}}
\left (
\int _{-\infty}^\infty
h(x_k, \mabfz _k-\mabffy _k(\mabfm _k))
^{p_k}\, dx_k
\right )
a_{k-1}(l_k,\mabfm _k)^{p_k}
\right )^{\frac 1{p_k}}
\\[1ex]
=
\left (
\sum _{\mabfm _{k}}
f_k(\mabfz _k-\mabffy _k(\mabfm _k))
^{p_k}
\left (
\sum _{l_k}a_{k-1}(l_k,\mabfm _k)^{p_k}
\right )
\right )^{\frac 1{p_k}}
\\[1ex]
=
\left (
\sum _{\mabfm _{k}}
f_k(\mabfz _k-\mabffy _k(\mabfm _k))
^{p_{0,k}}
a_{k}(\mabfm _k)^{p_{0,k}}
\right )^{\frac 1{p_{0,k}}},
\end{multline*}
and \eqref{Eq:ConvIndEst1} follows in this case
as well.

\par

It remain to consider the case $p_k\le p_{0,k-1}$
and $J_k=[0,1]$. Then $\fy _k(\mabfm _{k-1})=0$, and
by similar arguments as above we get
\begin{multline*}
g_k(\mabfz _k)
\\[1ex]
\lesssim
\left (
\int _{0}^1
\left (
\sum _{\mabfm _{k-1}}
h(x_k, \mabfz _k-\mabffy _k(\mabfm _k))
^{p}a_{k-1}(l_k,\mabfm _k)^{p}
\right )^{\frac {p_k}{p}}\, dx_k
\right )^{\frac 1{p_k}}
\\[1ex]
\lesssim
\left (
\int _{0}^1
\sum _{\mabfm _{k-1}}
\left (
h(x_k, \mabfz _k-\mabffy _k(\mabfm _k))
^{p}a_{k-1}(l_k,\mabfm _k)^{p}
\right )^{\frac {p_k}{p}}\, dx_k
\right )^{\frac 1{p_k}}
\end{multline*}
\begin{multline*}
=
\left (
\sum _{\mabfm _{k-1}}
\left (
\int _{0}^1
h(x_k, \mabfz _k-\mabffy _k(\mabfm _k))
^{p_k}\, dx_k
\right )
a_{k-1}(l_k,\mabfm _k)^{p_k}
\right )^{\frac 1{p_k}}
\\[1ex]
=
\left (
\sum _{\mabfm _{k}}
f_k(\mabfz _k-\mabffy _k(\mabfm _k))
^{p_{0,k}}
a_{k}(\mabfm _k)^{p_{0,k}}
\right )^{\frac 1{p_{0,k}}},
\end{multline*}
and \eqref{Eq:ConvIndEst1}, and thereby the result
follow.
\end{proof}

\par

\end{document}